\documentclass[12pt]{amsart}
\usepackage{a4wide}
\usepackage{amssymb}
\usepackage{mathrsfs}

\theoremstyle{plain}
\newtheorem{theorem}{Theorem} \newtheorem{lemma}{Lemma}[section]
\newtheorem{propo}{Proposition}[section]

\newtheorem{corol}{Corollary}[section]

\newcommand{\N}{\mathbb{N}}
 \newcommand{\Z}{\mathbb{Z}}

\newcommand{\G} {\Gamma} 
 
\newcommand{\Fix}{\mbox{Fix}}  \newcommand{\dist}{\mbox{dist}}
 
\newcommand{\Stab}{\mbox{Stab}}\newcommand{\Cl}{\mbox{Cl}}
\newcommand{\cR}{\mathcal{R}}
\newcommand{\cT}{\mathcal{T}} \newcommand{\cL}{\mathcal{L}}
\newcommand{\cB}{\mathcal{B}}
 \newcommand{\cN}{\mathcal{N}}

\newcommand{\Si}{S^0_\infty} \newcommand{\Ainf}{A^0_\infty}
\newcommand{\Sch}{\mbox{Sch}}
\newcommand{\Syml}{S^0_{[\,l,\infty]}} \newcommand{\Ainfl}{A^0_{[\,l,\infty]}}
\title[Invariant subsets of the space of subgroups]{Invariant subsets of the space of subgroups, equational compactness
and the weak equivalence of actions\footnote{
\,, Research partly sponsored by MTA Renyi ``Lendulet'' Groups
and Graphs Research Group}}

\author{G\'abor Elek and Konrad Kr\'olicki}

\begin{document}

\begin{abstract}
Equationally compact subgroups of countable groups were introduced
by Banaschewski. For all known cases the orbit closure of such a subgroup
is a countable subset
in the space of subgroups and has finite Cantor-Bendixson rank. We show that
there exists a finitely generated group $\G$ such that for any countable
ordinal $\alpha$ we have an equationally compact subgroup $H\subset \G$
for which the Cantor-Bendixson rank of the orbit closure of $H$ equals to 
$\alpha+2$. Then we give an explicite construction of continuum many
equationally compact subgroups of $\Gamma$ such that the
associated ergodic Bernoulli shift actions are pairwise weakly incomparable.
We also answer two questions on equational compactness posed 
by Prest and Rajani.

\end{abstract} \maketitle
\section{Introduction}
It is well-known that fields
are {\bf equationally compact}, that is, if we have a system $E_i, i\in I$ of 
equations over a 
field $K$
$$\sum_{j} a_{ij} x_j=m_i$$
and for any finite subset of $I$ the equations $E_i$ can be solved
simultaneously, then all the equations can be solved simultaneously. 
Now let $\G$ be a countable group acting on a set $X$ by permutations.
Let $\Fix(s)$ be the fixed point set of $s\in\Gamma$. 
Following Banaschewski \cite{Bana},
we
say that a $\G$-action is equationally compact, if we have a subset $S$ of $\G$
and for any finite subset $T$ of $S$, $\cap_{s\in T} \Fix(s)$ is non-empty, then
$\cap_{s\in S} \Fix(s)$ is non-empty.
A subgroup $H$ of $\G$ is equationally compact (or PIP, \cite{Prest}) 
if the left action on
$\G/H$ is equationally compact. The following proposition is quite straightforward
and is left for the reader.
\begin{propo} \label{propo1}
The subgroup $H\subset G$ is equationally compact if any of the following three conditions hold:
\begin{itemize}
\item $H$ is a finite extension of a normal subgroup (in particular, if $H$
is finite or normal).
\item The normalizator subgroup of $H$ has finite index in $G$.
\item $H$ is malnormal.
\end{itemize}
\end{propo}
\noindent
On the other hand, Banaschewski proved (\cite{Bana}, Proposition 6.) that
the free group of infinite generators has non-equationally compact subgroups.

\noindent
Let $\G$ be a countable group and $\{0,1\}^{\G}$ be the set of subsets of $\G$ 
with
the Tychonoff-topology. The set of all subgroups,
$S(\G)$ forms a closed, invariant (under the conjugate action) subspace
of $\{0,1\}^{\G}$. If $\G$ acts on a set $X$, then the set $\{\Stab(x): x\in X\}=M(\G,X)$
is an invariant subspace of $S(\G)$. The following proposition is easy to
prove.
\begin{propo} \label{easy}
Let $\G$ act on the set $X$ and let $M(\G,X)\subset S(\G)$ be the corresponding
invariant subspace of $S(\G)$. Then, the action is equationally compact
if and only if for any subgroup $K$ in the closure of $M(\G,X)$, there exists
an element of $L\in M(\G,X)$ such that $K\subseteq L$.
Particularly, $H\subset \G$ is equationally compact if for any $K$ in
the orbit closure of $H$, there exists a conjugate of $H$, $L=gHg^{-1}$ such
that $K\subseteq L$.
\end{propo}
The first goal of the paper is to answer two queries of Prest and Rajani 
\cite{Prest} concerning equational compactness by proving the following
two theorems.
\begin{theorem} \label{theorem1}
The only equationally compact subgroup of the finitary symmetric 
group $\Si$ on
$\N$ are the finite subgroups and the group of even permutations $A^0_\infty$.
\end{theorem}
\begin{theorem} \label{theorem2}
There exists a countable group $\G$ acting on a set $X$, such that
the action is equationally compact, but for any $x\in X$, $\Stab(x)$ is
not an equationally compact subgroup of $\G$.
\end{theorem}
\noindent
Let $H$ be a subgroup of a countable group $\G$. Then one can consider
the Cantor-Bendixson rank of its orbit closure (\cite {Cornul}). In some
sense the Cantor-Bendixson rank measures the complexity of a subgroup, how
far they are from being normal. 
Notice that the orbit closures of  all the
subgroups described in Proposition \ref{propo1} are countable sets and their
Cantor-Bendixson ranks are finite. 
One of the two main results of our paper is the following theorem.
\begin{theorem} \label{theorem3}
Let $\G=\Z_2 * \Z_2 *\
 Z_2$ be the free product of three cyclic groups.
Then, for any countable ordinal $\alpha$, there exists an equationally
compact subgroup $H$ of $\G$ such that orbit closure
of $H$ is countable and its Cantor-Bendixson rank is $\alpha+2$.
\end{theorem}
\noindent
Finally, we apply our techniques for the construction of weakly incomparable essentially free and ergodic generalized Bernoulli actions of
the group $\G_5=\Z_2 * \Z_2 * \Z_2 * \Z_2 * \Z_2$  (see Section \ref{weak}).
\begin{theorem}\label{theorem4}
There exist uncountable many equationally compact subgroups $\{H_\alpha\}_{\alpha\in I}$ of $\Gamma_5$ such that
the associated generalized Bernoulli actions $\G_5\curvearrowright ([0,1],\lambda)^{\Gamma_5/H_\alpha}$ are pairwise
weakly incomparable.
\end{theorem}

\section{Equationally compact subgroups of the finitary symmetric group}
Let $\Si=\cup^\infty_{n=1} S_n$ be the finitary symmetric group on the natural numbers.
That is, the group of permutations fixing all but finitely many elements.
The goal of this section is to prove Theorem \ref{theorem1} by showing that
the list of equationally compact subgroups of  $\Si$
contains only the set of finite groups and the alternating subgroup of
even permutations $\Ainf$.
Before getting into the proof let us fix some notations. Let  $\Syml\subset S^0_\infty$ be the subgroup of elements
fixing the set $\{1,2,\dots,l-1\}$. Let $\Ainfl=\Syml\cap A^0_\infty$. For a permutation $\gamma\in S^0_\infty$, we define
$s(\gamma)$ as the maximum of $k$'s for which $\gamma(k)\neq k.$ 
\begin{propo}
Let $H$ be an equationally compact subgroup of $S^0_\infty$. Then one of the following two conditions are satisfied.
\begin{enumerate}
\item There exists $l\geq 0$ such that $H\cap \Syml=\{e\}$.
\item There exists $l\geq 0$ such that $\Ainfl\subset H$.
\end{enumerate}
\end{propo}
\proof
Let $H\subset \Si$ be a subgroup such that neither of the two conditions above are satisfied.
Let $\kappa_1H\kappa^{-1}_1,  \kappa_2H\kappa^{-1}_2,\dots$ be an enumeration of the conjugates of $H$. Inductively,
we will pick elements
$\{\gamma_n\}^\infty_{n=1}\subset S^0_\infty, \{\delta_n\}^\infty_{n=1}\subset S^0_\infty$ such that
\begin{itemize}
\item $\{\delta_1,\delta_2,\dots,\delta_n\}\subset \gamma_n H \gamma_n^{-1}\,.$
\item $\delta_n\notin \kappa_n H \kappa_n^{-1}\,.$
\end{itemize}
\noindent
Hence the subgroup $H$ cannot be equationally compact.
Suppose that $\{\gamma_i\}^n_{i=1},  \{\delta_i\}^n_{i=1}$ has already been constructed and for any $1\leq i \leq n$
\begin{itemize}
\item $\{\delta_1,\delta_2,\dots,\delta_i\}\subset \gamma_i H \gamma_i^{-1}\,.$
\item $\delta_i\notin \kappa_i H \kappa_i^{-1}\,.$
\end{itemize}
\noindent
Let
$$l=\max(\max_{1\leq i \leq n} s(\gamma_i), \max_{1\leq i \leq n} s(\delta_i), \kappa_{n+1})+1\,.$$
\noindent
Since a conjugacy class always generates a normal subgroup, there exists a non-unit conjugacy class $C$ of $\Syml$
 such that $H\cap C$ is a proper subset of $C$.
Let $\delta_{n+1}\in C\backslash H, \rho_{n+1}\in H\cap C$.
Then, we have $\gamma\in \Syml$ such that
$\gamma\rho_{n+1}\gamma^{-1}=\delta_{n+1}$. Let $\gamma_{n+1}=\gamma\gamma_n.$
By the definition of $l$, we have that $\gamma$ commutes with $\{\gamma_i\}^n_{i=1}, \{\delta_i\}^n_{i=1}$
and $\kappa_{n+1}$, hence
\begin{itemize}
\item $\delta_{n+1}\notin \kappa_{n+1}H\kappa_{n+1}^{-1}\,.$
\item $\delta_i\in \gamma_{n+1} H \gamma_{n+1}^{-1}$, whenever $1\leq i \leq n\,.$
\end{itemize}
\noindent
Therefore, $H$ is not equationally compact. \qed
\begin{lemma}
If $H$ is equationally compact and contains $\Ainfl$ for some $l>0$, then either $H=\Si$ or $H=\Ainf$.
\end{lemma}
\proof
If $\Ainfl\subseteq H$, then for any $k\geq 1$ there
exists a conjugate of $H$, $\gamma H \gamma^{-1}$ such that the subgroup $A_k$ is contained in  $\gamma H \gamma^{-1}$.
Hence, if $H$ is equationally compact, then it must contain the whole group $\Ainf=\cup^\infty_{k=1} A_k$.
Therefore, $H=\Si$ or $H=\Ainf$. \qed

\vskip 0.2in
\noindent
The following proposition finishes to proof of Theorem \ref{theorem1}.
\begin{propo} If there exists $l\geq 1$ such that $H\cap \Syml=\{e\}$, then $H$ is finite.
\end{propo}
\proof Suppose that $H$ is an infinite subgroup of $\Si$.
\begin{lemma}
There exists an infinite subset $\{\gamma_n\}^\infty_{n=1}\subset H$ such that 
for any $n\geq 1$, $\gamma_n(1)=1\,.$
\end{lemma}
\proof First, let us suppose that there exists $k\geq 1$ and an infinite subset $\{\delta_n\}^\infty_{n=1}\subset H$
such that $\delta_n(1)=k.$
Let $\gamma_n=\delta_n^{-1}\delta_1.$ Then for  any $n\geq 1$, $\gamma_n(1)=1\,.$ 
If such $k$ does not exist, then we have an increasing sequence of positive integers $\{k_n\}^\infty_{n=1}$ and an
infinite subset $\{\gamma_n\}^\infty_{n=1}\subset H$ such that
\begin{itemize}
\item $\gamma_n(1)=k_n.$
\item $k_n>s(\gamma_i)$, whenever $1\leq i \leq n-1.$
\end{itemize}
\noindent
Then for any $n\geq 1$ and $1\leq i \leq n-1$, $\gamma_n^{-1}\gamma_i\gamma_n(1)=1$, hence our lemma follows. \qed
\vskip 0.2in
\noindent
Now let $s$ be a positive integer and
 $\{\delta_n\}^\infty_{n=1}\subset H$ be an infinite set of permutations such that $\delta_n(j)=j$ if $1\leq j\leq s$.
\begin{lemma}
There exists an infinite subset $\{\gamma_n\}^\infty_{n=1}\subset H$ such that
$\gamma_n(j)=j$, if $1\leq j\leq s+1\,.$
\end{lemma}
\proof
Again, if there exists $k\geq 1$ and an infinite subset $\{\rho_n\}^\infty_{n=1}\subset H$ such that
\begin{itemize}
\item $\rho_n(j)=j,$ if $1\leq j \leq s$.
\item $\rho_n(s+1)=k.$
\end{itemize}
\noindent
then the set $\{\gamma_n=\rho^{-1}_i\rho_1\}^\infty_{n=1}$ will satisfy the condition of our lemma.
On the other hand, if such $k$ does not exist then we have an increasing sequence of positive integers
 $\{k_n\}^\infty_{n=1}$ and an
infinite subset $\{\delta_n\}^\infty_{n=1}\subset H$ such that
\begin{itemize}
\item $\delta_n(j)=j$, if $1\leq j \leq s$.
\item $\delta_n(s+1)=k_n\,.$
\item $ k_n\geq s(\delta_i)$ if $1\leq i \leq n-1$.
\end{itemize}
\noindent
Hence $\delta^{-1}_n\delta_i\delta_n(j)=j$, if $1\leq  i \leq n$, $1\leq j \leq s+1 $. 
Thus our lemma follows. \qed
\vskip 0.2in
\noindent
By induction, we can construct infinitely many elements $\{\gamma_n\}^\infty_{n=1}\subset H$ such that $\gamma_n\in \Syml$,
in contradiction with the fact that $H\cap \Syml=\{e\}\,.$ \qed
\section{Tree subgroups and equational compactness} \label{treesub}
Let $T$ be a tree of vertex degrees at most three with edges
properly colored by the letters $a$, $b$ and $c$, that is,
adjacent edges are colored differently. From now on 
(until Section \ref{basiccon} all trees will be
considered to be properly $(a,b,c)$-edge-colored.
Trees are Schreier-graphs, the associated action of $\Gamma=\Z_2*\Z_2*\Z_2$
on the vertex set $V(T)$
is given the following way.
\begin{itemize}
\item If for $x\in V(T)$ there exists an edge $e(x,y)$ colored by $a$, then
$ax=y$.
\item If such edge does not exist, then $a(x)=x.$
\end{itemize}
\noindent
We define the action of the generators $b$ and $c$ in a similar fashion.
Let $w=p_n p_{n-1} p_{n-2}\dots p_1$, where $p_i=a,b$ or $c$ be a
reduced word, $x\in V(T)$. Then $w\in \Stab_T(x)$, if there exists
a closed walk $\{v_0, v_1,\dots, v_n\}$ in $T$ such that
\begin{itemize}
\item $v_0=v_n=x.$
\item $p_i(v_{i-1})=v_i\,.$
\end{itemize}
\noindent
If $x\in V(T)$, then $\Stab_T(x)$ is called a {\bf tree subgroup} of $\Gamma$. 
The orbit of $H=\Stab_T(x)$ in the space of
subgroups $S(\Gamma)$ is the set $\{\Stab(y)\,\mid\, y\in V(T)\}.$
Now we describe the orbit closure of $H$.
First, let us recall the notion of convergence and limits for rooted trees.
A rooted tree $(T,\rho)$ is an $(a,b,c)$-edge-colored tree with
a distinguished vertex $\rho\in V(T)$. The distance of two rooted trees is
defined the following way.
$$d_\tau\left( (T_1,\rho_1), (T_2,\rho_2)\right)=2^{-k}$$
\noindent
if the balls $B_k(T_1,\rho_1)$ and $B_k(T_2,\rho_2)$ are isomorphic as 
rooted, colored graphs, but  $B_{k+1}(T_1,\rho_1)$ and $B_{k+1}(T_2,\rho_2)$ 
are not isomorphic. The space $\cR\cT$ of rooted trees is a totally
disconnected compact space
with respect to the metric $d_\tau$. We will denote by $\cT$
the set of all $\{a,b,c\}$-trees up to isomorphism.
\begin{lemma}
The sequence $\{T_n,\rho_n\}^\infty_{n=1}$ converges to $\{T,\rho\}^\infty_{n=1}$
if and only if $\Stab_{T_n}(\rho_n)$ converges to $\Stab_T(\rho)$ in the
space of subgroups $S(\Gamma)$.
\end{lemma}
\proof
Let $\{T_n,\rho_n\}^\infty_{n=1}$ be a sequence of rooted trees
converging to $(T,\rho)$. Let $w\in \Gamma$ be a reduced word of length $k$.
By definition, if $n$ is large enough, then
$B_k(T_n,\rho_n)\cong B_k(T,\rho)$, hence $w(\rho)=\rho$ if and only if
$w(\rho_n)=\rho_n$ for large enough $n$. Therefore, $\Stab_{T}(\rho)$ is the 
limit of the sequence $\Stab_{T_n}(\rho_n)$ in the space  $S(\Gamma)$.
Conversely, let $H$ be the limit of the sequence
$\Stab_{T_n}(\rho_n)$ in $S(\Gamma)$. Let $W_k$ be the finite set of 
reduced words of length $k$ in $\Gamma$.
Then there exists an integer $N>0$ such that
if $N\leq n$, then an element $w\in W_n$ fixes $\rho_n$ if and only
if it fixes $\rho$.
Hence, the Schreier-graph $\Sch(\Gamma/H)$ must be a tree and
$\{T_n,\rho_n\}^\infty_{n=1}$ converges to the rooted tree
$\Sch(\Gamma/H, H)$ in the space of rooted trees $\cR\cT$.  \,\, \qed
\vskip 0.2in
\noindent
Before getting further, let us recall the notion of a {\bf branch}.
If $T$ is a tree and $x\in V(T)$, then $y,z\in V(T)\backslash\{ x\}$ is called
branch equivalent if the shortest paths connecting $y$ to $x$ resp.
$z$ to $x$ have a joint edge. For a fixed $x$ Tte equivalence classes induce subtrees of
$T$ and they are called branches. For any edge $e$ adjacent to $x$ there
exists exactly one branch containing $e$. The vertex $x$ is called the root of the branch. Now let $T_3$ be the
infinite tree of vertex degrees $3$ with its unique $(a,b,c)$-edge coloring.

\noindent
A {\bf decorated tree} $T^M_3$ is constructed the following way.
Let $M\subseteq E(T_3)$ be a set of edges of $T_3$ colored by the letter $a$.
For each chosen edge $e(x,y)\in M$ we will have two extra vertices $x_e,y_e$
such that $x$ is connected to $x_e$ by an edge colored by $a$,
$x_e$ is connected to $y_e$ by an edge colored by $b$ and $y_e$ is
connected to $y$ by an edge colored by $a$ (and then we delete the original edge $e(x,y))$. It is easy to see that the
resulting tree $T^M_3$ is properly $(a,b,c)$-colored,
the ``old'' vertices have valency $3$ and the ``new''  vertices have
valency $2$. The following lemma is trivial.
\begin{lemma}
For any decorated tree we have an $(a,b,c,D)$-coloring of the $3$-tree
 such that edges adjacent to an edge colored by $D$ can be colored
by either $b$ or $c$ (we call such edge-colorings {\bf good}). Conversely, if we have a good $(a,b,c,D)$-coloring
of the edges of the $3$-tree, then there is an associated decorated tree
$T^M_3$, where $M$ is the set of edges colored by $D$.
\end{lemma}
\noindent
It is easy to see that if
$\{(T_n,\rho_n)\}^\infty_{n=1}$ is a convergent sequence of good rooted
$(a,b,c,D)$-trees, then its limit $(T,\rho)$ is a good tree as well. Let $\cR \cT_D$ denotes the
compact space of rooted $(a,b,c,D)$-trees and $\cT_D$ denote the set of all $(a,b,c,D)$-trees. For a good 
$(a,b,c,D)$-tree $T$, $\phi(T)$ denotes the decorated $(a,b,c)$-colored tree associated to
$T$. 
The following lemma is easy to prove.
\begin{lemma}\label{coloring} The sequence
$\{(T_n,\rho_n)\}^\infty_{n=1}$ is converging to $(T,\rho)$ in $\cR \cT_D$ if and only if
$\{(\phi(T_n),\rho_n)\}^\infty_{n=1}$ is converging to $(\phi(T),\rho)$ in $\cR \cT$\,. 
\end{lemma}
\vskip 0.2in
\noindent
Now we define  {\bf sequential dominance} that is a crucial notion of our paper. 
Let $T$ and $S$ be good $(a,b,c,D)$-trees and $x\in V(T_3)$. We say that $S$ dominates $T$ via $x$ if there exists a 
branch $\cB$ in the underlying infinite tree $T_3$ such that
\begin{itemize}
\item
$x$ is not in $\cB$.
\item Outside the branch $\cB$ the colorings of $S$ and $T$ are identical.
\item There is no $D$-colored edge in $E(T)\cap \cB$.
\end{itemize}
\vskip 0.1in
\noindent
Now, suppose that $T$, $\{S_n\}^\infty_{n=1}$ are good $(a,b,c,D)$ trees, $x\in V(T)$ such that
$S_1$ dominates $T$ via $x$ and $S_n$ dominates $S_{n-1}$ via $x$ for any $n\geq 2$.
Then we have
$$E_D(T)\subset E_D(S_1)\subset E_D(S_2)\subset\dots\,$$
\noindent
where $E_D(T)$ is the set of $D$-colored edges in $T$. So we can consider the good $(a,b,c,D)$-tree $S$, where
$E_D(S)=\cup^\infty_{n=1} E_D(S_n)$. In this case, we say that $S$ sequentially dominates $T$ via $x$.
\begin{lemma}\label{workhorse}
Let $(R,x)$ and $(S,y)$ be good $(a,b,c,D)$-trees. Suppose that there exists a good $(a,b,c,D)$-tree $T$
such that
\begin{itemize}
\item $S$ sequentially dominates $T$ via $y$.
\item $(R,x)$ and $(T,y)$ are isomorphic as rooted, colored trees.
\end{itemize}
\noindent
Then $\Stab_{\phi(R)}(x)\subset \Stab_{\phi(S)}(y)\,.$ If the conditions of this lemma are satisfied then we say that the
rooted tree $(S,y)$ sequentially dominates $(R,x)$. Note that for simplicity we write $y$ instead of $\phi(y)$
since the vertex $y$ can be clearly identified in $\phi(S)$.
\end{lemma}
\proof
It is enough to show that $\Stab_{\phi(T)}(y)\subset \Stab_{\phi(S_1)}(y)$.
Let $w$ be a reduced word that fixes $y$ in $\phi(T)$. Since all the vertex degrees of $\phi(T)$ in the branch $\cB$
are $3$, any closed walk in $\pi(T)$ that goes through $y$ and has an edge in $\cB$ must have a turning point. Hence,
the walk induced by $w$ does not have an edge in $\cB$. Therefore, $w\in \Stab_{\phi(S_1)}(y)\,.$
Inductively,
$$\Stab_{\phi(T)}(y) \subset \Stab_{\phi(S_1)}(y) \subset \Stab_{\phi(S_2)}(y) \subset \cdots$$
\noindent
Since $\Stab_{\phi(S)}(y)=\cup^\infty_{n=1}  \Stab_{\phi(S_n)}(y)$, our lemma follows. \qed
\vskip 0.2in
\noindent
The following proposition is our main tool for constructing equationally compact subsgroups in $\Gamma$.
\begin{propo} \label{main}
Let $S$ be an $(a,b,c,D)$-tree. Let
$\{(T_\alpha,x_\alpha)\}_{\alpha\in I}$ be the set of all elements of the closure of the set $\{S,y\}_{y\in S}$ in $\cR \cT_D$.
Suppose that for any $\alpha\in I$, there exists $x_\alpha\in V(T_\alpha)$ and $y_\alpha \in V(S)$ such that $(S,y_\alpha)$
sequentially dominates $(T_\alpha,x_\alpha)$. Then for each $y\in V(S)$, $\Stab_{\phi(S)}(y)$ is an equationally compact
subgroup of $\Gamma$.
\end{propo}
\proof By Proposition \ref{easy}, we only need to prove
that if $\{(\phi(S),y_i)\}^\infty_{i=1}$ converges to $(Q,z)$ in $\cR \cT$, then
$\Stab_Q(z)\subset \Stab_{\phi(S)}(y)$ for some $y\in V(S)$.
We can suppose that the vertices $\{y_i\}$ has degree $3$. Indeed, if the degrees (for large $i$) are $2$, then
we can substitute $y_i$ by an adjacent vertex $y'_i$ having vertex degree $3$. Then $\{(\phi(S),y'_i)\}^\infty_{i=1}$ will
converge to $(Q,z')$, where $z'$ is adjacent to $z$. Now, if $\Stab_Q(z')\subset\Stab_{\phi(S)}(y')$ for some
$y'\in V(S)$, then $\Stab_Q(z)\subset \Stab_{\phi(S)}(y)$ for some $y\in V(S)$.
By Lemma \ref{coloring}, if $\{(\phi(S),y_i)\}^\infty_{i=1}$ converges to $(Q,z)$ in $\cR \cT$, then
$\{(S,y_i)\}^\infty_{i=1}$ converges to $(T,z)$ in $\cR \cT_D$ and $\phi(T)=Q$.
By the condition of our proposition, there exists $z"\in V(T)$ such that the rooted tree $(T,z")$ is
sequentially dominated by $(S,y')$ for some $y'\in V(S)$.
Hence, by Proposition \ref{workhorse}, $\Stab_{\phi(T)}(z")\subset \Stab_{\phi(S)}(y')$. That is,
$$\Stab_Q(z)=\Stab_{\phi(T)}(z)\subset \Stab_{\phi(S)}(y)\,\,\mbox{ for some}\,\,  y\in V(S)\quad\qed$$
\vskip 0.2in
\noindent
The following example is intended to illustrate the use of Proposition \ref{main}.
We call a good $(a,b,c,D)$-tree $S$ {\bf sparse}, if for any $n\geq 1$, there exists only
finitely many pairs of edges $E_D(S)$ having distance at most $n$.
\begin{propo} \label{sparse}
Let $S$ be a sparse tree. Then for any vertex $y\in\phi(S)$, $\Stab_{\phi(S)}(y)$ is equationally compact.
\end{propo}
\proof Suppose that the sequence $\{(S,x_i)\}^{\infty}_{i=1}$ converges to $(T,x)$ in $\cR \cT_D$. Recall, that it means that for any
$r\geq 1$, the ball $B_r(T,x)$ is rooted-colored-isomorphic to $B_r(S,x_i)$ if $i$ is large enough. We have three cases.
\vskip 0.1in
\noindent
{\bf Case 1.}
The sequence  $\{x_i\}^\infty_{i=1}$ is bounded. Then $(T,x)$ is isomorphic to $(S,y)$ for some vertex $y\in V(S)$.
\vskip 0.1in
\noindent
{\bf Case 2.}\, There exists $r\geq 1$ such that
$B_r(S,x_i)$ contains exactly one $D$-colored edge. 
Then the limit tree $(T,x)$ has exactly one $D$-colored edge. Clearly, $(T,x)$ is sequentially dominated by $(S,y)$ for some $y\in S$.
\vskip 0.1in
\noindent
{\bf Case 3.}\, For all $r\geq 1$, the balls $B_r(S,x_i)$ does not contain $D$-colored edges if $i$ is large enough. Then $(T,x)$ has
no $D$-colored edges, so $(T,x)$ is sequentially dominated by $(S,y)$ for all $y\in V(S)$ and $\Stab_{\phi(T)}(x)=\{e\}$.

\vskip 0.1in
\noindent
Hence, the conditions of Proposition \ref{main} are satisfied. Therefore, 
for any vertex $y\in\phi(S)$, $\Stab_{\phi(S)}(y)$ is equationally compact. \qed
\section{The proof of Theorem \ref{theorem3}}.
Before getting into the proof of Theorem \ref{theorem3} let us introduce some definitions to make the construction of good $(a,b,c,D)$-
trees easier.
A bi-infinite $(a,b)$-path is a graph on the vertex set $\{x_n\}_{n\in \Z}$ such that $x_n$ and $x_m$ is adjacent if and only if
$|n-m|=1$ and the edge $e(x_{2k},x_{2k+1})$ is colored by $a$, the edge $e(x_{2k},x_{2k-1})$ is colored by $b$.
An infinite $(c,b)$-path is a graph on the vertex set $\{y_n\}_{n\geq 0}$ such that  $y_n$ and $y_m$ is adjacent if and only if
$|n-m|=1$ and the edge $e(y_{2k},y_{2k+1})$ is colored by $c$, the edge $e(y_{2k},y_{2k-1})$ is colored by $b$.
Now we  build our $(a,b,c,D)$-tree. We start with the standard infinite $(a,b,c)$-colored tree $T_3$. 
\vskip 0.2in
\noindent
{\bf Step 1.} We fix a bi-infinite 
$(a,b)$-subpath in $T_3$ on the vertices $\{x_n\}_{n\in \Z}$. 
\vskip 0.2in
\noindent
{\bf Step 2.} Let us consider the infinite $(c,b)$-path starting from $x_0$. The vertices of the path are  $\{y^0_n\}^\infty_{n=0}$,
where $y^0_0=x_0$.
\vskip 0.2in
\noindent
{\bf Step 3.} Now for any $k>0$ we consider the infinite  $(c,b)$-path starting from $x_{2^k}$. The vertices of this path are
$\{y^k_n\}^\infty_{n=0}$, where $y^k_0=x_{2^k}$.

\vskip 0.2in
\noindent
{\bf Step 4.} For any $i\geq 1$, we choose a subset $L_i$ of the positive integers. If $s\in L_i$, then we recolor the $a$-colored
edge adjacent to the vertex $y^i_{2^s}$ by the color $D$.
\vskip 0.2in
\noindent
{\bf Step 5.} We recolor the $a$-colored edges adjacent to the vertices $\{y^0_{2^k}\}^\infty_{k=1}$ by $D$.
\begin{propo} \label{ltree} For all choices of the family $\cL$ of the sets $\{L_i\}^\infty_{i=1}$ the resulting tree $S_\cL$ satisfies the conditions of Proposition \ref{main},
hence $\Stab_{\phi(S)}(z)$ is equationally compact if $z\in V(S)$.
\end{propo}
\proof
Again, suppose that the
sequence $\{(S,z_i)\}^\infty_{i=1}$ converges to $(T,q)$ in $\cR \cT_D$.
We need to show that $(T,q)$ is sequentially dominated by $(S,y)$ for some
$y\in V(S)$. The three possible cases described in Proposition \ref{sparse}
are already handled. The underlying limit trees are in these cases $S, \hat{T}_3$ and $T_3$, where 
$\hat{T}_3$ is the tree with one single $D$-colored edge. Before considering the fourth case, let us define
a specific tree. Let $L\in \{0,1\}^{\N}$ be representing a subset of the positive integers. 
Let $p=p_0,p_1,p_2,\dots$
be a $(c,b)$-path in $T_3$ starting from the vertex $p$. Then recolor by $D$
the $a$-colored edge adjacent to $p_i$, if $i=2^j$, where $j\in L$. The resulting
tree is denoted by $T_L$.
\vskip 0.1in
\noindent
{\bf Case 4.} For some $r\geq 1$, the balls $B_r(S,z_i)$ contains at least two $D$-colored edges provided
that $i$ is large enough.

\noindent
Hence, there exists some integer $d\geq 0$ such that $\dist_S(z_i,x_{n_i})=d$ for some increasing 
sequence $\{n_i\}^\infty_{i=1}$, where $\dist_S$ is the shortest path distance in $S_{\cL}$. First,
let us suppose that $z_i=x_{n_i}$ that is $d=0$.
The following lemma is an easy consequence of the definition of rooted tree convergence.
\begin{lemma} \label{tihonoff}
The sequence $\{L_{n_i}\}^\infty_{i=1}\subset \{0,1\}^{\N}$ is convergent in the Tychonoff-topology of pointwise
convergence and $(T,q)\cong (T_L,p)$, where $L$ is the limit of  $\{L_{n_i}\}^\infty_{i=1}$.
\end{lemma}
\noindent
By our construction, $(T_L,p)$ is sequentially contained in $(S,x_0)$. Now, if $d\neq 0$, then 
$(T,x)\cong (T_L,p')$, where $p'$ is some vertex of $T_L$. Hence, $(T,x)$ is still sequentially
contained by $(S,x_0)$. This finishes the proof of Proposition \ref{ltree}. \qed
\vskip 0.2in
\noindent
Now we are ready to finish the proof of Theorem \ref{theorem3}.
Let $\cN=\{N_j\}^\infty_{j=1}\subset \{0,1\}^\N$ be a countable closed set with Cantor-Bendixson
rank $\alpha$. It is well-known that such set exists. Choose the sets $\{L_i\}^\infty_{i=1}$ in such a way
that each $L_i$ equals to some $N_j$ and for 
each $j$ we have infinitely many $i$
so that $L_i$ equals to $N_j$. We denote by $\Cl(S_\cL)$ the closure of $\{S_\cL,y\}_{y\in S_\cL}$ in
$\cR \cT_D$.
The following proposition finishes the proof of Theorem \ref{theorem3}.
\begin{propo}
The Cantor-Bendixson rank of  $\Cl(S_\cL)$ is $\alpha+2$ if $\alpha$ is an infinite ordinal and
$\alpha+3$ if $\alpha$ is a finite ordinal.
\end{propo}
\proof
As we have see in the proof of  Proposition \ref{ltree}, $\Cl(S_\cL)$ is the union of the families
$\{(S_\cL,y)\}_{y\in S_\cL}$, $\cup^\infty_{j=1} \{(T_{N_j},z)\}_{z\in T_{N_j}}$, $\{(\hat{T}_3,x)\}_{x\in \hat{T}_3}$
and the element $(T_3,p)$. Let $\{M_k\}^\infty_{k=1}\subset \{0,1\}^\N$ and 
$Q:=\cup^\infty_{k=1}\{(T_{M_k},z)\}_{z\in T_{M_k}}.$ Again, by Proposition \ref{ltree}, $\Cl(Q)$ is the
union of the families 
$\cup_{L\in\overline {\cup^\infty_{k=1} M_k}} \{(T_{L},y)\}_{y\in T_L}, \{(\hat{T}_3,x)\}_{x\in \hat{T}_3}$ and the
element  $(T_3,p)$. Hence, the isolated points of  $\Cl(S_\cL)$ denoted by $\Cl_0(S_\cL)$ are exactly the
elements $\{(S_\cL,y)\}_{y\in \cL}$. The isolated points of $\Cl(S_\cL)\backslash \Cl_0(S_\cL)$ denoted by $\Cl_1(S_\cL)$
are $\cup_{N_j\in \cN_0} \{(T_{N_j},z)\}_{z\in T_{N_j}}$, where $\cN_0$ is the set of isolated points in $\cN$.
For an ordinal $\gamma$, define $\cN_\gamma$ to be the isolated points of 
$\cN\backslash \cup_{\beta<\gamma} \cN_\beta$. Similarly,
define $\Cl_{\gamma}$ to be the isolated points of $\Cl(S_\cL)\backslash \cup_{\beta<\gamma}\Cl_\beta(S_\cL) $.
By transfinite induction, we can see that for finite ordinals $\beta$
$$\Cl_{\beta+1}(S_\cL)=\cup_{N_j\in\cN_\beta} \{(T_{N_j},z)\}_{z\in T_{N_j}}\,,$$
\noindent
and for infinite ordinals $\beta$
$$\Cl_{\beta}(S_\cL)=\cup_{N_j\in\cN_\beta} \{(T_{N_j},z)\}_{z\in T_{N_j}}\,.$$
Therefore by the definition of the set $\cN$, if $\alpha$ is a finite ordinal
$$\Cl(S_\cL)\backslash \Cl_{\alpha+1}(S_\cL)=\{(\hat{T}_3,x)\}_{x\in \hat{T}_3} \cup (T_3,p)\,,$$
\noindent 
and if $\alpha$ is an infinite ordinal
$$\Cl(S_\cL)\backslash \Cl_{\alpha}(S_\cL)=\{(\hat{T}_3,x)\}_{x\in \hat{T}_3} \cup (T_3,p)\,.$$
Hence for a finite ordinal $\Cl_{\alpha+3}$ consists of one single element $(T_3,p)$.
For an infinite ordinal  $\Cl_{\alpha+2}$ consists of one single element $(T_3,p)$. \qed
\section{Minimal systems}
The goal of this section is to prove Theorem \ref{theorem2}.
First, let us recall the notion of minimal subshifts of the Bernoulli shift space.
Let $\{a,b,c\}^\Z$ be the set of all $\{a,b,c\}$-valued functions $\sigma$ on the integers with
the natural $\Z$-action
$$t_n(\sigma)(a)=\sigma(a-n)\,.$$
\noindent
A minimal subshift is a closed, invariant subspace $\Sigma\subset \{a,b,c\}^\Z$ such that the orbit closure
of any $\sigma\in\Sigma$ is $\Sigma$ itself.
Let $w=(q_k,q_{k-1},\dots,q_1)\in \Gamma=\Z_2*\Z_2*\Z_2$ and $\sigma\in  \{a,b,c\}^\Z$. We say that $n\in\Z$ {\it sees} $w$ if
$$\sigma(n-i)=q_i\quad\mbox{for any $1\leq i \leq k$.} $$
\noindent
It is easy to see (e.g. \cite{monod}) that the orbit closure of $\sigma$ is a minimal subshift if for any $w\in\Gamma$ that is seen by
some integer $n$, there exists $m_w>0$ such that the longest interval in $\Z$ without elements that see $w$ is shorter than $m_w$.
We call such a $\sigma$ a minimal sequence. A {\it good} minimal sequence is a minimal sequence that
does not contain the same letter consecutively. It is well-known that good minimal sequences exist for which
the associated subshift has the cardinality of the continuum. 
 Now, let us consider the line graph $L$ on $\Z$.
That is, $a,b\in \Z$ is connected if and only if $|a-b|=1$. Let $\sigma$ be a good  minimal sequence. 
Color the edge $(n,n+1)$ of $L$ by $\sigma(n)$.
Then we obtain the $\{a,b,c\}$-tree $L_\sigma$.
Let $\Sigma$ be the orbit closure of  $\sigma$.
By Lemma \ref{coloring} and the definition of a  good minimal sequence one can see immediately that if $\tau\in\Sigma$ and $n\in\Z$, then
$\Stab_{L_\tau}(n)$ is in the orbit closure of $\Stab_{L_\sigma}(0)$ in the space of subgroups $S(\Gamma)$. Conversely, any element
of the orbit closure of $\Stab_{L_\sigma}(0)$ is in the form of $\Stab_{L_\tau}(n)$, for some $\tau$ and $n$.
That is, 
$$\{\Stab_{L_\tau}(n)\}_{\tau\in\Sigma, n\in\Z}$$
\noindent
is a minimal $\Gamma$-system in $S(\Gamma)$. These systems are called a uniformly recurrent subgroup (URS) in \cite{GW}
(we thank Mikl\'os Ab\'ert to call our attention to this paper). Since there are continuum many minimal subshifts
in $\{a,b,c\}^\Z$, in this way we obtain continuum many URS's in $S(\Gamma)$ (see Theorem 5.1 \cite{GW}). 
Note that if $\Sigma$ and $\Sigma'$ represents the same URS, then either $\Sigma=\Sigma'$ or $\Sigma'=\Sigma^{-1}$, where
$\sigma^{-1}(n):=\sigma(-n)$, for any $\sigma\in\Sigma$. Clearly, an URS is an equationally compact set in $S(\Gamma)$, therefore
Theorem \ref{theorem2} follows from the lemma below.
\begin{lemma}
Let $\sigma,\tau$ be good elements of a minimal subshift $\Sigma$. Suppose that $\tau$ is neither a $\Z$-translate of $\sigma$ nor a
 $\Z$-translate of $\sigma^{-1}$. Then for any $n\geq 1$, $\Stab_{L_\sigma}(0)\not\subset \Stab_{L_\tau}(n)$.
\end{lemma}
\proof
Let $n\geq 1$ and let $w=(q_k,q_{k-1},\dots,q_1)$ be the longest word such that
one of the following two conditions hold.
\begin{enumerate}
\item $\sigma(-i)=\tau(n-i)$ for any $1\leq i \leq k$. 
\item $\sigma(-i)=\tau(n+i)$ for any $1\leq i \leq k$. 
\end{enumerate}
\noindent
Without a loss of generality we can suppose that the first condition holds, $\sigma(-k)=c$, $\sigma(-k-1)=a$, $\tau(n-k-1)=b$.
Then clearly, $w^{-1}bw\in \Stab_{L_\sigma}(0)$. On the other hand, $w^{-1}bw\notin \Stab_{L_\tau}(n)$. Indeed, $bw(n)=n-k-1$, hence 
$w^{-1}bw(n)<n\,.$ \qed

\section{Weak equivalence of  actions I. (the basic construction)}
\label{basiccon}

In this section we construct
a continuum of equationally compact subgroups of the group $\Gamma_5=\Z_2*\Z_2*\Z_2*\Z_2*\Z_2.$
Let $\{a,b,c,d,e\}$ be free generators of order two in $\Gamma_5$.
Let $T_5$ be the infinite tree of vertex degrees five properly edge-colored by the symbols
$\{a,b,c,d,e\}$.
We will consider good $\{a,b,c,d,e,D\}$-colorings of $T_5$ and the associated $\Gamma_5$-actions,
where a $D$-colored edge is modified exactly the same way as in Section \ref{treesub},
by adding two vertices $p_x,p_y$ such that the edges
$(x,p_x)$ and $(p_y,y)$ are colored by $a$ and the edge $(p_x,p_y)$ is colored by $b$.
So, in the associated $\{a,b,c,d,e\}$-tree the degrees of the new vertices are two, the
degrees of the old vertices are five.

\vskip 0.2in
\noindent
{\bf The master subtree R}.  In order to build our master subtree $R$ in $T_5$ we need some
definitions.
A $(2^n,b,a)$-path in $T_5$ is a path
$(x_0,x_1,x_2,\dots,x_{2^n})$ such that all edges $(x_{2l},x_{2l+1})$ are colored by $b$ and all
edges $(x_{2l+1},x_{2l+2})$ are colored by $a$. Hence the path starts with a $b$-colored edge and ends
with an $a$-colored edge.
We define the $(2^n,c,a)$-paths similarly. They start with a $c$-colored edge and end with an $a$-colored edge.
Now, fix a vertex $p\in V(T_5)$. Pick the $(2,b,a)$-path and the $(2,c,a)$-path starting from $p$
to obtain two new endpoints. Then for each of the two endpoints pick the $(4,b,a)$-path and the $(4,c,a)$-path
starting from them to obtain four new endpoints. In the $n$-th step we see $2^{n-1}$-endpoints and for each
such endpoints $x$ we pick the $(2^n,b,a)$-path and the $(2^n,c,a)$-path starting from $x$ to obtain $2^n$ new endpoints.
Inductively, we build the infinite subtree $R$, as the union of all the chosen paths above.

\vskip 0.2in
\noindent
{\bf The codes}. A code $C$ is an infinite sequence $\{a_n\}^\infty_{n=1}$ such that if $n\geq 1$, then $a_i=b$ or $a_i=c$.
Using our master tree $R$, for each code $C$ we build a good $(a,b,c,d,e,D)$-tree $T_C$ and the associated $(a,b,c,d,e)$-tree
$S_C$. Note that $T_C$ will be recoloring of $T_5$ and not a recoloring of the master subtree $R$.
\vskip 0.1in
\noindent
{\bf Step 1.} We recolor by $D$ the last edge of the $(2,a_1,a)$-path starting from $p=p_C$ we denote by $x^C_1$ the endpoint of the
{\it other} $2$-path (a $(2,b,a)$-path or a $(2,c,a)$-path) starting from $p$.
\vskip 0.1in
\noindent
{\bf Step 2.} We recolor by $D$ the last edge of the $(4,a_2,a)$-path starting from $x^C_1$. Then we denote the endpoint of the other $4$-path by
$x^C_2$.
\vskip 0.1in
\noindent
{\bf Step n.} We recolor by $D$ the last edge of the $(2^n,a_n,a)$-path starting from $x^C_{n-1}$ and 
denote by $x_n^C$ the endpoint of the
other $2^n$-path.
\vskip 0.1in
\noindent
Inductively, we construct the good $(a,b,c,d,e,D)$-tree $T_C$. \\ Since the associated $(a,b,c,d,e)$-tree $S_C$ is 
sparse, by Proposition \ref{sparse}, $\Gamma_C:=\Stab_{S_C}(p)$ is an equationally
compact subgroup of $\Gamma_5$. The following technical proposition is crucial for the proof of Theorem \ref{theorem4}.
\begin{propo}\label{hresz}
Let $C=(a_1,a_2,\dots)$ and $C'=(a'_1,a'_2,\dots)$ be different codes. Then there exists a finitely generated 
non-amenable subgroup $H\subset \Gamma_C$ such
that no finite subset of $V(S_{C'})$ is invariant under the action of $H$.
\end{propo}
\proof
We prove our propostion using two lemmas.
\begin{lemma}
$$\Stab_{S_C}(p_C)\not\subset \Stab_{S_{C'}}(p_{C'})$$
\end{lemma}
\proof
Since $C\neq C'$, there exist $q\in V(S_C)$, $q'\in V(S_{C'})$ and $n\geq 1$  such that
\begin{itemize}
\item The paths $(p_C,q)$ and $(p_{C'},q')$ are colored-isomorphic paths of length $n$.
\item The degree of $q$ is two and the degree of $q'$ is five.
\end{itemize}
\noindent
Let $w\in\Gamma_5$ be the word of length $n$ such that
$w(p_C)=q$ and $w(p_{C'})=q'$. Then $w^{-1}e w\in \Stab_{S_C}(p_C)$ and $w^{-1}e w\notin \Stab_{S_{C'}}(p_{C'}).$\,\,\,\qed
\vskip 0.2in
\noindent
Let $\gamma=a_1aeaa_1$, $\delta=\overline{a_1}a a_2 a a_2 a e a a_2 a a_2 a \overline{a_1}$, where 
$\overline{a_1}=b$ if $a_1=c$ and $\overline{a_1}=c$ if $a_1=b$. Obviously, $\gamma,\delta\in \Stab_{S_C}(p_C)$.
\begin{lemma}
Let $F$ be a finite subset of $V(S_{C'})$ containing at least one element that does not equal to $p_{C'}$. Then
$F$ is not invariant under the group generated by $\gamma$ and $\delta$.
\end{lemma}
\proof
Suppose that $q\in F$ and $q$ is not a vertex of degree $2$ or a vertex from the master subtree $R$. Then
either the $a_1$-branch or the $a_2$-branch of $q$ contains only vertices of degree five. Hence,
either the set $\{(\gamma\delta)^n(q)\}^\infty_{n=1}$ or the set  $\{(\delta\gamma)^n(q)\}^\infty_{n=1}$ is infinite.
So, we can suppose that all the elements of $F$ are either vertices of degree $2$ or they are from the master subtree.
Let $s$ be one of the furthest elements from $p_{C'}$ in $F$. Then either $\gamma(s)$ or $\delta(s)$ are further
from $p_{C'}$ than $s$. Hence $F$ cannot be invariant under the group generated by $\gamma$ and $\delta$.\,\,\,\qed
\vskip 0.2in
\noindent
Now, let $H\subset \Stab_{S_C}(p_C)$ be generated by the set
$$\{ \gamma, \delta, a_1adaa_1,  a_1aeaa_1,  a_1afaa_1, \alpha\,\}\,,$$
\noindent
where $\alpha\in \Stab_{S_C}(p_C)\backslash \Stab_{S_{C'}}(p_{C'})\,.$ Then $H$ is finitely generated, non-amenable subgroup
of $\Stab_{S_C}(p_C)$, without non-empty invariant subsets in $V(S_{C'})$. \,\,\,\qed
\vskip 0.2in
\noindent
Now we recall the notion of amenable actions \cite{Tsankov}. An action $\alpha:\Gamma \curvearrowright X$ of a countable group $\Gamma$
on a countable set $X$ is called {\bf amenable} if there exists a sequence of finite subsets in $X$, the so-called F{\o}lner sets 
$\{F_n\}^\infty_{n=1}$ such that for any $g\in\Gamma$
$$\lim_{n\to\infty} \frac{|\alpha(g)(F_n)\cup F_n |}{| F_n|}=1\,.$$
\noindent
Suppose that $\Gamma$ is generated by the finite set $\{g_1,g_2,\dots,g_r\}$. Then the action $\alpha$ is non-amenable if and only if
there exists an $\epsilon>0$ such that for any finite set $F\subset X$
$$|\cup^r_{j=1} \alpha(g_j)(F)\cup F|\geq (1+\epsilon)|F|\,.$$
\noindent
The left action of a group $\Gamma$  on itself or on an overgroup of $\Gamma$  is amenable if and only if the group $\Gamma$ is amenable. 
\begin{propo} \label{technical}
Let $C\neq C'$ be two codes and $H\subset \Stab_{S_C}(p_C)$ be the finitely generated non-amenable subgroup of $\Gamma_5$ as in
Proposition \ref{hresz}. then the restricted action of $H$ on the set $V(S_{C'})$ is non-amenable.
\end{propo}
\proof We start with a technical lemma.
\begin{lemma}
Let $\{G_n\}^\infty_{n=1}$ be a sequence of connected, induced, finite subgraphs in $V(S_{C'})$. Then
$$\liminf_{n\to\infty}\frac{|\cup^r_{j=1} \beta(h_j)(V(G_n))\cup V(G_n)|}{ |V(G_n)|}>1\,.$$
\noindent
Here, $\beta$ denotes the restricted action of $H$.
\end{lemma}
\proof
For any $n\geq 1$, we consider the finite induced subgraph $H_n\subset T_5$ constructed the 
following way.
If $(p,q,r,s)$ is a path in $S_{C'}$ such that the degrees of $q$ and $p$ are two in $S_{C'}$ with respect to the standard generating
system $\{a,b,c,d,e\}$ (that is, this path was substituted for a $D$-colored edge) and the path is in $G_n$, then let $(p,s)\in E(H_n)$.
If only $(p,q)$ or $(r,s)$ is an edge of $G_n$, then $(p,s)\notin E(H_n)$.  On the other hand, if
the degrees of $p$ and $q$ are five in $S_{C'}$, then $(p,q)\in H_n$ if and only if $(p,q)\in G_n$. So we substitute new paths with
a single edge and cut off hanging edges starting from new vertices.  Hence, there is a bijection between
$V_5(G_n)$ and $V(H_n)$, where 
$$V_5(G_n)=\{x\in V(G_n)\,\mid\, \deg_{S_{C'}}(x)=5\}\,.$$
\noindent
Since $S_{C'}$ is sparse (that is for any $n\geq 1$ there are only finitely many pairs of degree $2$ vertices of distance less than
$n$),\begin{equation} \label{e1}
\lim_{n\to\infty}\frac{|\cup^r_{j=1} \alpha(h_j)(V(H_n))\cup V(H_n)|}{|\cup^r_{j=1} \beta(h_j)(V(G_n))
\cup V(G_n)|}=1\,,
\end{equation}
\noindent
where $\alpha$ denotes the left action of $H$ on the overgroup $\Gamma_5$. By the non-amenability of $\alpha$, 
$$\lim_{n\to\infty}\frac{|\cup^r_{j=1} \alpha(h_j)(V(H_n))\cup V(H_n)|}{|V(H_n)|}>1\,,$$
hence by (\ref{e1}) our proposition follows. \qed
\begin{corol}
There exist $\delta>0$ and $C>0$ such that
for any finite, connected induced subgraph $L\subset G$ such that $|V(L)|\geq C$
$$|\cup^r_{j=1} \beta(h_j)(V(L))\cup V(L)|>(1+\delta) |V(L)|\,.$$ \end{corol}
\noindent
Now, we finish the proof of the proposition. Let $M\subset V(S_{C'})$ be a finite subset and $L_M$ be the subgraph induced by $M$.
Let $\{L_q\}^s_{q=1}$ be the connected components of $L_M$. By Proposition \ref{hresz} $V(L_j)$ is not invariant under
$H$, hence for any $1\leq q\leq s$
$$|\cup^r_{j=1} \beta(h_j)(V(L_q))\cup V(L_q)|\geq \frac{C+1}{C} |V(L_q)|\,,$$
provided that $|V(L_q)|\leq C$.
Hence we have the following inequality. For any $1\leq q \leq s$,
\begin{equation} \label{ucso}
|\cup^r_{j=1} \beta(h_j)(V(L_q))\backslash V(L_q)|\geq m |V(L_q)|\,,
\end{equation}
\noindent
where $m=\min(\frac{1}{C},\delta)$.

Observe that for any vertex $x\in L_M$, there exist at most $r$ values of $q$ such that
$$x\in \cup^r_{j=1} \beta(h_j)(V(L_q))\backslash V(L_q)\,,$$
hence by (\ref{ucso}), we have that
$$|\cup^r_{j=1} \beta(h_j)(V(L_M))\cup V(L_M)|\geq
(1+\frac{m}{r}) |V(L_M)|\,.\quad\qed$$

\section{Weak equivalence of actions II. (the proof of Theorem \ref{theorem4})} \label{weak}

The notion of weak equivalence was introduced by Kechris \cite{Kechris} and
since then it has been studied extensively. Let $\alpha:\Gamma \curvearrowright (M,\lambda), 
\beta:\Gamma \curvearrowright (M,\lambda)$ be measure preserving actions of a countable group $\Gamma$ on a 
standard probability space
$(M,\lambda)$. We say that $\alpha$ weakly contains $\beta$, $\alpha \succeq \beta$ if for any
finite measurable partition of $M$, $M=\cup^n_{i=1}B_i$, a finite subset $\{g_j\}^m_{j=1}\subset \Gamma$ and a constant $\epsilon>0$,
there exists a partition $M=\cup^n_{i=1}A_i$ such that for any $1\leq j \leq m$, $1\leq k,l \leq n$
$$\left| \lambda(\alpha(g_j)(A_k)\cap A_l)- \lambda(\beta(g_j)(B_k)\cap B_l) \right| \leq \epsilon\,.$$
\noindent
The actions $\alpha,\beta$ are {\bf weakly incomparable} if $\alpha \not\succeq \beta$ and $\beta \not\succeq \alpha$.
 It was proved in \cite{Abert} that
for certain groups $\Gamma$ there exist uncountably many pairwise weakly incomparable free and ergodic probability measure 
preserving (p.m.p.) actions of $\Gamma$.
We shall prove in Theorem \ref{theorem4} that for $\Gamma_5$ there exist uncountably many generalized shifts associated to
subgroups of $\Gamma_5$ that are pairwise weakly incomparable. We will use a result of Kechris and Tsankov [Theorem 1.2]\cite{Tsankov} in a crucial way.
Let $\alpha:\Gamma \curvearrowright (M,\lambda)$ be a p.m.p.action. A sequence of measurable subsets in $M$, $\{A_n\}^\infty_{n=1}$ are called
asymptotically invariant with respect to the action $\alpha$ if for any $g\in \Gamma$
$$\lim_{n\to \infty} \lambda\left(\alpha(g)(A_n)\backslash A_n\right)=0\,.$$
The action $\alpha$ is called {\bf strongly ergodic} if for any asymptotically invariant sequence $\{A_n\}^\infty_{n=1}$ 
$$\lim_{n\to \infty} \lambda(A_n)\lambda(M\backslash A_n)=0\,.$$
\noindent
Clearly, a strongly ergodic action of $\Gamma$ cannot weakly contain a non-ergodic action of $\Gamma$.  Now, let
$\beta:\Gamma \curvearrowright  X$ be an action of $\Gamma$ on a countable set by permutations.
The associated generalized Bernoulli action $\hat{\beta}$ is defined
on the product space $\prod_{x\in X} ([0,1],\nu)$, where $\nu$ is the Lebesgue measure. The action is defined by
$$\hat{\beta}(g)(F)(x)=F(\beta(g^{-1})(x))\,,$$
where $F\in [0,1]^X, g\in \Gamma$.
It is easy to see that $\hat{\beta}$ is essentially free if for any $e\neq g\in \Gamma$, there exist infinitely many $x\in X$
such that $\beta(g)(x)\neq x\,.$
According to the theorem of Kechris and Tsankov, the associated generalized Bernoulli action 
$\hat{\beta}:\Gamma \curvearrowright \prod_{x\in X} ([0,1],\nu)$ is strongly ergodic if and only if the action $\beta$ is non-amenable.
Now we can finish the proof of Theorem \ref{theorem4}. 
Let $\prod^\infty_{n=1} \{b,c\}$ be the set of all codes. For each code $C$ we constructed an $\{a,b,c,d,e\}$-colored tree $S_C$ and the
associated $\Gamma_5$-action $\alpha_C:\Gamma_5\curvearrowright V(S_C)$. Let $C\neq C'$ be two codes and $H\subset \Gamma_5$ be the subgroup of
$\Stab_{S_C}(p_C)$ as in Proposition  \ref{hresz}. Then the associated generalized Bernoulli action restricted to $H$, 
$$\hat{\alpha}_C\mid_H: H \curvearrowright \prod_{v\in V(S_C)}([0,1],\nu)$$
\noindent
is clearly non-ergodic.  On the other hand, by Proposition \ref{technical}
and the aforementioned theorem of Kechris and Tsankov,
$$\hat{\alpha}_{C'}\mid_H: H \curvearrowright \prod_{v\in V(S_{C'})}([0,1],\nu)$$
\noindent
is strongly ergodic, therefore $\hat{\alpha}_{C'}\not\succeq \hat{\alpha}_C$.
Hence Theorem \ref{theorem4} follows. \qed

\vskip 0.4in
\noindent
G\'abor Elek, Lancaster University \\
g.elek@lancaster.ac.uk \\
\vskip 0.2in
\noindent
K. Kr\'olicki, Lancaster University \\
k.krolicki@lancaster.ac.uk
\end{document}